\documentclass[reqno]{amsart}%
\usepackage{float}
\usepackage{amsmath}
\usepackage{graphicx}
\usepackage{latexsym}
\usepackage{amsfonts}
\usepackage{amssymb}%
\theoremstyle{plain}
\newtheorem{theorem}{Theorem}
\newtheorem{corollary}[theorem]{Corollary}
\newtheorem{lemma}[theorem]{Lemma}
\newtheorem{proposition}[theorem]{Proposition}
\theoremstyle{definition}

\newtheorem{definition}[theorem]{Definition}

\newtheorem{remark}[theorem]{Remark}
\newtheorem*{remark*}{Remark}

\begin{document}
\title[LD for sums on a GW process]{Large deviations for sums\vspace{10pt}\\defined on a Galton-Watson process\vspace{40pt}}
\author[Fleischmann]{Klaus Fleischmann}
\address{Weierstrass Institute for Applied Analysis and Stochastics, Mohrenstr.\ 39, \\
D--10117 Berlin, Germany}
\email{fleischm@wias-berlin.de}
\urladdr{http://www.wias-berlin.de/$\sim$fleischm}
\author[Wachtel]{Vitali Wachtel}
\address{Weierstrass Institute for Applied Analysis and Stochastics, Mohrenstr.\ 39, \\
D--10117 Berlin, Germany}
\email{vakhtel@wias-berlin.de}
\urladdr{http://www.wias-berlin.de/$\sim$vakhtel}
\thanks{WIAS preprint No.\ 1135 of May 23, 2006,\quad ISSN 0946\thinspace
--\thinspace8633,\quad LD26.tex}
\thanks{Supported by the DFG}
\thanks{Corresponding author: Klaus Fleischmann}
\subjclass{Primary 60\thinspace J\thinspace80; Secondary 60\thinspace F\thinspace10}
\maketitle

\thispagestyle{empty}
\setcounter{page}{0}\newpage

$\qquad$\bigskip

\begin{quotation}
\noindent\textsc{Abstract.} In this paper we study the large deviation
behavior of sums of i.i.d.\ random variables $X_{i}$ defined on a
supercritical Galton-Watson process $Z.$ We assume the finiteness of the
moments $\,EX_{1}^{2}$\thinspace\ and $\,EZ_{1}\log Z_{1\,}.$\thinspace\ The
underlying interplay of the partial sums of the\thinspace\ $X_{i}$ and the
lower deviation probabilities of $Z$ is clarified. Here we heavily use lower
deviation probability results on $Z$ we recently published in
\cite{FleischmannWachtel2006.LowerDev.AIHP}.

{\small
\tableofcontents
}
\end{quotation}

\section{Introduction and results}

\subsection{Motivation\label{SS.motive}}

Let $\,Z=(Z_{n})_{n\geq0}$\thinspace\ denote a \emph{Galton-Watson process}\/
with offspring law $\,\{p_{k};\,k\geq0\}$. We will assume that $Z$%
\thinspace\ is \emph{supercritical}: $\,m:=\sum_{k=1}^{\infty}kp_{k}%
\in(1,\infty).$\thinspace\ As a rule we start with $\,Z_{0}=1.$

A basic task in statistical inference of Galton-Watson processes is the
estimation of the offspring mean $m.$ Let us recall at this place the
well-known Lotka-Nagaev estimator $\,Z_{n+1}/Z_{n}$\thinspace\ of $\,m$ due to
A.V.\ Nagaev \cite{NagaevAV1967}. If $\,\sigma:=(\mathbf{V}\!\mathrm{ar}%
Z_{1})^{1/2}\in(0,\infty),$ then for every $\,x\in\mathbb{R}$,
\begin{equation}
\lim_{n\uparrow\infty}\mathbf{P}\Bigl(m^{n/2}\bigl(\tfrac{Z_{n+1}}{Z_{n}%
}-m\bigr)<x;\ Z_{n}>0\Bigr)\ =\ \int_{0}^{\infty}\Phi\bigl(\tfrac{xu^{1/2}%
}{\sigma}\bigr)\,w(u)\,\mathrm{d}u, \label{Fluct1}%
\end{equation}
where $w$ denotes the density function of the a.s.\ limit variable
$W:=\lim_{n\uparrow\infty}m^{-n}Z_{n}$ restricted to $(0,\infty),$%
\thinspace\ and $\Phi$ is the standard normal distribution function,
\begin{equation}
\Phi(y)\,:=\,\frac{1}{\sqrt{2\pi}}\int_{-\infty}^{y}\mathrm{e}^{-z^{2}%
/2}\,\mathrm{d}z,\qquad y\in\mathbb{R}.
\end{equation}

The study of the ratio $Z_{n+1}/Z_{n}$ has attracted the attention of several
researchers in recent years, since it can also be used for estimating
important parameters such as the amplification rate and the initial size in a
quantitative polymerase chain reaction experiment; see Jacob and Peccoud
\cite{JacobPeccoud1996,JacobPeccoud1998}.

Fix $\,k\geq0.$\thinspace\ In a finer description of the Galton-Watson model,
let $Z_{n}(k)$ denote the number of particles in the $n^{\mathrm{th}}$
generation having exactly $k$ children. Then, on the event $\,\{Z_{n}%
>0\},$\thinspace\ results for the estimator $\,\tilde{p}_{k}(n):=Z_{n}%
(k)/Z_{n}$\thinspace\ of $\,p_{k\,},$\thinspace\ which hold analogously to
(\ref{Fluct1}), had been provided by Pakes \cite[Theorems~5 and 6]{Pakes1975}.

The mentioned results from \cite{NagaevAV1967} and \cite{Pakes1975} can be
seen from a unified point of view as follows. Independently of $\,Z,$%
\thinspace\ let $\,X=(X_{n})_{n\geq1}\,\ $denote a family of
\emph{i.i.d.\ random variables with mean zero and variance in} $(0,\infty
).$\ Let $\,n\geq0.$\thinspace\ Put $\,S_{n}:=X_{1}+\ldots+X_{n\,}.$%
\thinspace\ On the event $\,\{Z_{n}>0\},$\thinspace\ the random variable $\,$%
\begin{equation}
R_{n}\,:=\,S_{Z_{n}}/Z_{n} \label{not.Rn}%
\end{equation}
is well-defined. For convenience, \emph{we agree}\/ that an event involving
$\,R_{n}$\thinspace\ is always tacitly assumed to be included in
$\,\{Z_{n}>0\}.$\thinspace\ For instance, $\,\mathbf{P}(R_{n}<x)$%
\thinspace\ means $\,\mathbf{P}\bigl(R_{n}<x;\ Z_{n}>0\bigr)$\thinspace\ more
carefully written. If now $X_{1}$ coincides in law with $Z_{1}-m,$ then, for
$\,n$\thinspace\ fixed, $\,R_{n}$\thinspace\ coincides in law with
$Z_{n+1}/Z_{n}-m$\thinspace\ on the event $\,\{Z_{n}>0\}.$\thinspace\ On the
other hand, if $X_{1}$ takes on the value $1-p_{k}$ with probability $p_{k}$
(for $\,k$\thinspace\ fixed) and $-p_{k}$ otherwise, then $R_{n}=\tilde{p}%
_{k}(n)-p_{k}$\thinspace\ in law on the event $\,\{Z_{n}>0\}$\thinspace\ for
fixed $\,n\,.$

Sums such as $S_{Z_{n}}$ arise also in models of polymerase chain reactions
with mutations, see Piau \cite{Piau2004}.\smallskip

From now on, we work with the more general meaning of $\,R_{n\,},$%
\thinspace\ based on $\,(X,Z),$\thinspace\ as introduced in (\ref{not.Rn}).
Clearly, we have the following \emph{strong law of large numbers:}%
\begin{equation}
R_{n}\rightarrow0\quad\text{a.s. \thinspace as}\quad n\uparrow\infty.
\end{equation}
Moreover, using methods from \cite{NagaevAV1967} and \cite{Pakes1975}, one can
easily verify the following \textquotedblleft\emph{normal deviation
probabilities}\textquotedblright\ for $R_{n}:$
\begin{equation}
\lim_{n\uparrow\infty}\mathbf{P}\bigl(m^{n/2}R_{n}<x\bigr)\ =\ \int
_{0}^{\infty}\Phi\bigl(\tfrac{xu^{1/2}}{\sigma}\bigr)\,w(u)\,\mathrm{d}%
u,\qquad x\in\mathbb{R}, \label{Fluct1'}%
\end{equation}
where $\,\sigma:=(\mathbf{E}X_{1}^{2})^{1/2}$\thinspace\ from now on. Let
$\,\varepsilon_{n}>0.$\thinspace\ In the case $\,\varepsilon_{n}%
m^{n/2}\rightarrow\infty,$\thinspace\ this implies the following simple large
deviation probabilities for $\,R_{n}:$%
\begin{equation}
\lim_{n\uparrow\infty}\mathbf{P}(R_{n}\geq\varepsilon_{n})=0.
\end{equation}
But the main task of large deviation theory is to determine the \emph{rate}\/
of this convergence. Clearly, one of the reasons to be interested in large
deviation probabilities comes from statistical applications. Firstly, these
probabilities describe the quality (error probabilities) of many tests. On the
other hand, a question concerning the Bahadur efficiency of estimators leads
also to the large deviation problem.

For the particular model $X_{1}\overset{\mathcal{L}}{=}Z_{1}-m,$%
\thinspace\ the special case $\varepsilon_{n}\equiv\varepsilon$ is more or
less studied in the literature. In fact, Athreya \cite{Athreya1994} proved
that if $\,p_{0}=0,$ $\,p_{1}>0,$\thinspace\ and $\,\mathbf{E}Z_{1}%
^{2\alpha+\delta}<\infty$\thinspace\ for some $\,\delta>0$, where $\alpha
\in(0,\infty)$ denotes the so-called Schr\"{o}der constant [see
(\ref{not.gamma.alpha}) below], then
\begin{equation}
\lim_{n\uparrow\infty}\,m^{\alpha n}\,\mathbf{P}\!\left(  _{\!_{\!_{\,}}%
}|R_{n}|\geq\varepsilon\right)  \quad\text{exists finitely.} \label{Athreya1}%
\end{equation}
On the other hand, using asymptotic properties of harmonic moments of
$Z_{n\,}$, Ney and Vidyashankar \cite{NeyVidyashankar2003} found the rate of
$\,\mathbf{P}\!\left(  _{\!_{\!_{\,}}}|R_{n}|\geq\varepsilon\right)  $ under
the weaker assumption that $\,\mathbf{P}(Z_{1}\geq j)\sim aj^{1-\eta}%
$\thinspace\ as $j\uparrow\infty$, for some $\,\eta>2$\thinspace\ and $a>0.$
The same authors proved in \cite{NeyVidyashankar2004} a version of a large
deviation principle for $R_{n}$ conditioned on $Z_{n}\geq v_{n}$ with numbers
$v_{n}\rightarrow\infty$; see also Rouault \cite{Rouault2000}.\smallskip

The \emph{purpose of the present paper}\/ is to study the rate of convergence
of large deviation probabilities of $\,R_{n}\geq\varepsilon_{n}$\thinspace\ in
the more interesting case $\,\varepsilon_{n}\rightarrow0$\thinspace\ as
$n\uparrow\infty$ (working with our more general setting of $\,R_{n}%
).$\thinspace\ For this we heavily relay on results on lower deviation
probabilities of $\,Z,$ we recently derived in
\cite{FleischmannWachtel2006.LowerDev.AIHP}. In the next subsection we briefly
recall what we need from that paper.

Note that large deviation probabilities in the case $\varepsilon
_{n}\rightarrow0$ are needed, for instance, for testing two close hypotheses,
i.e.\ when the distance between the hypotheses tends to zero as the size of
the sample gets larger and larger.

\subsection{Lower deviation probabilities for $Z$\label{SS.lower}}

We start with recalling the following basic notation, reflecting a crucial
dichotomy for supercritical Galton-Watson processes.

\begin{definition}
[\textbf{Schr\"{o}der and B\"{o}ttcher case}]\label{D.dichot}For our
supercritical offspring di\-stribution we distinguish between the
\emph{Schr\"{o}der}\/ and the \emph{B\"{o}ttcher}\/ case, in dependence on
whether\/ $\,p_{0}+p_{1}>0$\thinspace\ or $\,=\,0,$\thinspace
\ respectively.\hfill$\Diamond$
\end{definition}

Write $f$ $\,$for the generating function of our supercritical offspring law:
$\,f(s)=\sum_{j\geq0}$\thinspace$p_{j}s^{j},$\thinspace\ $0\leq s\leq
1.$\thinspace\ Let $q$\thinspace\ denote the extinction probability of
$\,Z,$\thinspace\
\begin{equation}
\text{set }\,\gamma\,:=\,f^{\prime}(q),\quad\text{and define}\,\ \alpha
\,\ \text{by}\quad\gamma\,=\,m^{-\alpha}. \label{not.gamma.alpha}%
\end{equation}
Note that $\,\gamma\in\lbrack0,1)$\thinspace\ and $\,\alpha\in(0,\infty
].$\thinspace\ Obviously, we are in the Schr\"{o}der case if and only if
$\,\gamma>0,$\thinspace\ if and only if $\,\alpha<\infty.$\thinspace\ In this
case, $\,\alpha$\thinspace\ is said to be the \emph{Schr\"{o}der constant.} We
also need the following notion.

\begin{definition}
[\textbf{Type }$(d,\mu)$]\label{D.type}We say the offspring distribution
\emph{is of type} $(d,\mu)$, if $\,d\geq1$\thinspace\ is the greatest common
divisor of the set $\,\{j-\ell:\,j\neq\ell,\ p_{j}p_{\ell}>0\},$%
\thinspace\ and $\,\mu\geq0$\thinspace\ is the minimal \thinspace$j$%
\thinspace\ for which \thinspace$p_{j}>0$.\hfill$\Diamond$
\end{definition}

In the present paper, $\,(d,\mu)$\thinspace\ \emph{always refers}\/ to the
type of our offspring law. Recall that in the B\"{o}ttcher case, $\,\mu
=\min\{k\geq0:p_{k}>0\}\geq2.$\thinspace\ Here the \emph{B\"{o}ttcher
constant}\/ \thinspace$\beta\in(0,1)$ is defined by $\,\mu=m^{\beta}%
$.\thinspace\ 

We also \emph{always assume}\/ that the moment $\,\mathbf{E}Z_{1}\log Z_{1}%
$\thinspace\ is finite. Under this moment condition, the results of
\cite{FleischmannWachtel2006.LowerDev.AIHP} can be specialized to the
following two propositions.

\begin{proposition}
[\textbf{Schr\"{o}der case}]\label{P.Schroeder}In the Schr\"{o}der case, for
$\,k_{n}\leq m^{n}$ $\,$satisfying $\,k_{n}\rightarrow\infty,$\thinspace\ we
have
\begin{equation}
\sup_{k\in\lbrack k_{n},m^{n}]\ \text{with\ }k\equiv\mu\,(\mathrm{mod}%
\,d)}\,\left\vert \,\frac{m^{n}}{d\,w\bigl(k/m^{n}\bigr)}\,\mathbf{P}%
(Z_{n}=k)\ -\ 1\,\right\vert \;\underset{n\uparrow\infty}{\longrightarrow}\;0
\label{44}%
\end{equation}
and
\begin{equation}
\sup_{k\in\lbrack k_{n},m^{n}]}\left\vert \,\frac{\mathbf{P}(0<Z_{n}\leq
k)}{\mathbf{P}\bigl(0<W<k/m^{n}\bigr)}\ -\ 1\,\right\vert \;\underset
{n\uparrow\infty}{\longrightarrow}\;0.
\end{equation}
{}
\end{proposition}

\begin{proposition}
[\textbf{B\"{o}ttcher case}]\label{P.Boettcher}Suppose the B\"{o}ttcher case.
Then there exist positive constants $B_{1}$ and $B_{2}$ such that for all\/
$\,k_{n}\geq\mu^{n}$\thinspace\ with $\,k_{n}=o(m^{n})$\thinspace\ as
$\,n\uparrow\infty,$
\begin{subequations}
\label{10}%
\begin{align}
-B_{1}\  &  \leq\ \liminf_{n\uparrow\infty}\,(k_{n}/m^{n})^{\beta/(1-\beta
)}\log\mathbf{P}(Z_{n}\leq k_{n})\\
\  &  \leq\ \limsup_{n\uparrow\infty}(k_{n}/m^{n})^{\beta/(1-\beta)}%
\log\mathbf{P}(Z_{n}\leq k_{n})\,\leq\ -B_{2\,}.
\end{align}
The inequalities\/ \emph{(\ref{10}) }remain true if $\,\mathbf{P}(Z_{n}\leq
k_{n})$ is replaced by $\,m^{n}\,\mathbf{P}(Z_{n}=k_{n}),$\thinspace\ provided
that $\,k_{n}\equiv\mu^{n}\,(\mathrm{mod}\,d)$.
\end{subequations}
\end{proposition}

In order to explain the influence of lower deviation probabilities of
$\,Z_{n}$ on $\,R_{n}=S_{Z_{n}}/Z_{n\,},$\thinspace\ we look at the
decomposition,
\begin{equation}
\mathbf{P}(R_{n}\geq\varepsilon_{n})\ =\ \sum_{k=1}^{\infty}\mathbf{P}%
(Z_{n}=k)\,\mathbf{P}(S_{k}\geq\varepsilon_{n}k). \label{TotalProb}%
\end{equation}
Thus, in order to find the asymptotics of $\,\mathbf{P}(R_{n}\geq
\varepsilon_{n}),$\thinspace\ we need to determine the range of values of
$\,k$, which give the main contribution in decomposition (\ref{TotalProb}). As
we will see, this depends on parameters of the offspring law (as $\alpha$, for
instance) and, on the other hand, on the tail behavior of $X_{1\,}$. Here we
mention several possibilities. If $k$ is of order $m^{n}$ (the regime of
normal deviations for $Z_{n}$) and $\varepsilon_{n}^{2}m^{n}\rightarrow\infty
$, then $\varepsilon_{n}k$ is in the domain of large deviations of $S_{k}$. On
the other hand, if $k$ is of order $\varepsilon_{n}^{-2}$ (regime of normal
deviations for $S_{k}$), then $k$ is in the domain of lower deviations for
$Z_{n\,}.$\thinspace\ And finally, if $k/m^{n}\rightarrow0$ and $\varepsilon
_{n}^{2}k\rightarrow\infty$, then simultaneously we have lower deviations for
$Z_{n}$ and large deviations for $S_{k\,}.$

\subsection{Large deviations in the Schr\"{o}der case}

Recall that we always assume $\,\mathbf{E}Z_{1}\log Z_{1}<\infty$%
\thinspace\ and $\,\mathbf{E}X_{1}^{2}<\infty.$\thinspace\ As usual, we set
$\,X_{1}^{+}:=X_{1}\vee0.$

\begin{theorem}
[\textbf{Schr\"{o}der under a }$(1+\alpha)$\textbf{-moment condition on
}$X_{1}$]\label{T_DDev}\hfill Suppose \newline the Schr\"{o}der case and that
\begin{equation}
\mathbf{E(}X_{1}^{+})^{1+\alpha}\,<\,\infty\label{MomCond}%
\end{equation}
\emph{[}with $\,\alpha\in(0,\infty)$\thinspace\ the Schr\"{o}der constant
from\/ \emph{(\ref{not.gamma.alpha})]. }Moreover, assume that $\,\varepsilon
_{n}\rightarrow0$\thinspace\ and $\,\varepsilon_{n}^{2}m^{n}\rightarrow\infty
$\thinspace\ as $\,n\uparrow\infty$. Then
\begin{subequations}
\label{13}%
\begin{align}
0\,<\,V_{\ast\,}\Gamma_{\alpha}  &  \,\leq\,\liminf_{n\uparrow\infty
}\,\,\varepsilon_{n}^{2\alpha}\,m^{\alpha n}\,\mathbf{P}(R_{n}\geq
\varepsilon_{n})\label{LB}\\
\,  &  \,\leq\,\limsup_{n\uparrow\infty}\,\varepsilon_{n}^{2\alpha}\,m^{\alpha
n}\,\mathbf{P}(R_{n}\geq\varepsilon_{n})\,\leq\,V^{\ast}\Gamma_{\alpha
}\,<\,\infty, \label{UB}%
\end{align}
where
\end{subequations}
\begin{equation}
V_{\ast}\,:=\,\liminf_{u\downarrow0}\,u^{1-\alpha}w(u),\quad V^{\ast
}\,:=\,\limsup_{u\downarrow0}\,u^{1-\alpha}w(u) \label{not.V.star}%
\end{equation}
and
\begin{equation}
\Gamma_{\alpha}\ :=\ \frac{2^{\alpha-1}\,\Gamma(\alpha+1/2)}{\alpha\sqrt{\pi}%
}\ \sigma^{2\alpha}. \label{notC.alpha}%
\end{equation}
{}
\end{theorem}

\noindent Of course, here $\,\Gamma$\thinspace\ refers to the Gamma function.

Next we recall some known facts on the asymptotic behavior of supercritical
Galton-Watson processes in the Schr\"{o}der case. With $q$ and $\gamma
$\thinspace\ introduced in the beginning of Subsection~\ref{SS.lower} and with
$\,f_{n}$\thinspace\ denoting the iterates of $\,f,$\thinspace\ the following
limit exists:
\begin{equation}
\lim_{n\uparrow\infty}\,\frac{f_{n}(s)-q}{\gamma^{n}}\ =:\ \mathsf{S}%
(s)\ =:\ \sum_{j=0}^{\infty}\nu_{j}s^{j},\qquad0\leq s<1. \label{SConv'}%
\end{equation}
Hence,
\begin{equation}
\lim_{n\uparrow\infty}\gamma^{-n}\,\mathbf{P}(Z_{n}=k)\,=\,\nu_{k\,},\qquad
k\geq1. \label{SConv}%
\end{equation}
The Schr\"{o}der constant $\alpha<\infty$ describes the behavior of the
density function $\,w(u)$ as $u\downarrow0$. In fact, according to Biggins and
Bingham \cite{BigginsBingham1993}, there is a continuous, positive
multiplicatively periodic function $V$ such that
\begin{equation}
u^{1-\alpha}w(u)\,=\,V(u)+o(1)\quad\text{as }\,u\downarrow0. \label{w-asymp}%
\end{equation}
The function $V$ in (\ref{w-asymp}) can be replaced by a (positive) constant
$V_{0}$ if and only if
\begin{equation}
\mathsf{S}\!\left(  _{\!_{\!_{\,}}}\varphi(h)\right)  =\,V_{0}h^{-\alpha
},\qquad h\geq0, \label{V-condition}%
\end{equation}
where $\,\varphi$\thinspace\ denotes the Laplace transform of the limit random
variable $\,W$\thinspace\ (cf.\ Asmussen and Hering \cite[p.\thinspace
96]{AsmussenHering1983}. In this case, $\,V^{\ast}=V_{\ast}$\thinspace\ in
Theorem~\ref{T_DDev}, and we get the following conclusion.

\begin{corollary}
[\textbf{Schr\"{o}der under an additional regularity of }$Z$]\label{C1}
If\/\emph{ (\ref{V-condition})} holds, then under the assumptions of
Theorem\/\emph{~\ref{T_DDev}},
\begin{equation}
\lim_{n\uparrow\infty}\varepsilon_{n}^{2\alpha}\,m^{\alpha n}\,\mathbf{P}%
(R_{n}\geq\varepsilon_{n})\ =\ V_{0}\,\Gamma_{\alpha} \label{Cor}%
\end{equation}
\emph{[}with $\,\Gamma_{\alpha}$ from \emph{(\ref{notC.alpha})]}.
\end{corollary}

Under the assumptions of Theorem~\ref{T_DDev}, the sum at the right hand side
of (\ref{TotalProb}) is determined by those values of $k$ which are of order
$\varepsilon_{n}^{-2}.$ As we already mentioned, this corresponds to lower
deviations of $\,Z$ and normal deviations of $\,S_{k\,}.$\thinspace\ But what
happens if moment condition (\ref{MomCond}) fails? We are able to answer this
question under some regularity of the tail probabilities of $\,X_{1\,}%
.$\thinspace\ For this purpose, we say that $\,X_{1}$\thinspace\ \emph{has a
tail of index} $\,\theta,$\thinspace\ if for some constant $\,a>0,$%
\begin{equation}
\mathbf{P}(X_{1}\geq x)\,\sim\,a\,x^{-\theta}\quad\text{as }\,x\uparrow\infty.
\label{TailAssumption}%
\end{equation}
(Here the involved constant is always denoted by $\,a.)$

\begin{theorem}
[\textbf{Schr\"{o}der under heavier tails concerning }$X_{1}$]\label{T_LDev}%
Suppose that $\,1$ $<\alpha<\infty$\thinspace\ and that $\,X_{1}$%
\thinspace\ has a tail of index $\,\theta\in(2,1+\alpha).$\thinspace\ Define
$\,\varkappa:=(1+\alpha-\theta)/(2\alpha-\theta).$

\begin{enumerate}
\item[(a)] If\/ $\,\varepsilon_{n}m^{\varkappa n}\rightarrow0,$\thinspace\ but
$\,\varepsilon_{n}^{2}m^{n}\rightarrow\infty,$\thinspace\ then statements\/
\emph{(\ref{13})} remain valid.\vspace{2pt}

\item[(b)] If $\,\varepsilon_{n}m^{\varkappa n}\rightarrow\infty,$%
\thinspace\ then
\begin{equation}
\lim_{n\uparrow\infty}\varepsilon_{n}^{\theta}\,m^{(\theta-1)n}\,\mathbf{P}%
(R_{n}\geq\varepsilon_{n})\ =\ a\,I_{\theta}, \label{LargeDev}%
\end{equation}
where
\begin{equation}
I_{\theta}\ :=\ \frac{1}{\Gamma(\theta-1)}\int_{0}^{\infty}\varphi
(v)\,v^{\theta-2}\,\mathrm{d}v. \label{notC.theta}%
\end{equation}

\item[(c)] If\/ $\,\varepsilon_{n}m^{\varkappa n}\rightarrow\tau^{-1}%
\in(0,\infty),$\thinspace\ then%
\begin{align}
\hspace*{30pt}\tau &  ^{2\alpha}V_{\ast}\Gamma_{\alpha}+\tau^{\theta
}a\,I_{\theta}\ \leq\ \liminf_{n\uparrow\infty}\,m^{\alpha(\theta
-2)n/(2\alpha-\theta)}\,\mathbf{P}(R_{n}\geq\varepsilon_{n})\\
&  \leq\ \limsup_{n\uparrow\infty}m^{\alpha(\theta-2)n/(2\alpha-\theta
)}\,\mathbf{P}(R_{n}\geq\varepsilon_{n})\ \leq\ \tau^{2\alpha}V^{\ast}%
\Gamma_{\alpha}+\tau^{\theta}a\,I_{\theta}\nonumber
\end{align}
\emph{[}with $\,V_{\ast\,},V^{\ast}$ from\/ \emph{(\ref{not.V.star})} and
$\,\Gamma_{\alpha}$\thinspace\ from\/ \emph{(\ref{notC.alpha})].}
\end{enumerate}
\end{theorem}

Large deviations as in part (b) have a structure, different from that in
Theorem~\ref{T_DDev}. Here the main contribution comes from normal deviations
of $\,Z_{n}$ and large deviations of $\,S_{n\,}.$\thinspace\ In part (c) we
have a combination of regimes appearing in (a) and (b).

\begin{remark}
[\textbf{Critical value of }$\,\theta$]Theorems \ref{T_DDev} and \ref{T_LDev}
leave open the case that $\,X_{1}$\thinspace\ has a tail of index
$\,\theta=\alpha+1.$\thinspace\ Our methods allow to prove that (under
$\,1<\alpha<\infty)$\thinspace\ part (a) holds, if $\,\varepsilon
_{n}n^{1/(\alpha-1)}\rightarrow0.$\thinspace\ On the other hand, if
$\,\varepsilon_{n}n^{1/(\alpha-1)}\rightarrow\infty,$\thinspace\ then%
\begin{equation}
\lim_{n\uparrow\infty}n^{-1}\,\varepsilon_{n}^{1+\alpha}\,m^{\alpha
n}\,\mathbf{P}(R_{n}\geq\varepsilon_{n})\ =\ a\,J_{\alpha}%
\end{equation}
where%
\begin{equation}
J_{\alpha}\,:=\,\frac{1}{\Gamma(\alpha)}\int_{1}^{m}\mathsf{S}\!\left(
_{\!_{\!_{\,}}}\varphi(v)\right)  v^{\alpha-1}\,\mathrm{d}v.
\end{equation}
Finally, if $\,\varepsilon_{n}n^{1/(\alpha-1)}\rightarrow\tau^{-1}\in
(0,\infty),$\thinspace\ then a similar statement as in (c) is true.\hfill
$\Diamond$
\end{remark}

\subsection{Large deviations in the B\"{o}ttcher case\label{SS.Boe}}

As well-known, in the B\"{o}ttcher case the following limit
\begin{equation}
\lim_{n\uparrow\infty}\bigl(f_{n}(s)\bigr)^{(\mu^{-n})}\ =:\ \mathsf{B}%
(s),\qquad0\leq s\leq1, \label{BConv}%
\end{equation}
exists, is positive and continuous. From this it follows that in general
$f_{n}(s)$ does not converge as $n\uparrow\infty.$ But taking logarithms, we
have%
\begin{equation}
\lim_{n\uparrow\infty}\,\mu^{-n}\log f_{n}(s)\,=\,\log\mathsf{B}(s).
\end{equation}
On the other hand, our result on lower deviations in the B\"{o}ttcher case
(Proposition~\ref{P.Boettcher}) is also only for log-scaled probabilities.
These two facts explain the use of a logarithmic scaling in our following theorem.

\begin{theorem}
[\textbf{B\"{o}ttcher under light tails concerning }$X_{1}$]\label{Bottcher}%
\hfill Assume the B\"{o}tt-\newline cher case, that $\,\mathbf{E}%
\mathrm{e}^{h|X_{1}|}$ is finite for some $\,h>0,$\thinspace\ and that
$\,\varepsilon_{n}\rightarrow0$\thinspace\ as well as $\,\varepsilon_{n}%
^{2}m^{n}\rightarrow\infty$\thinspace\ as $\,n\uparrow\infty$. Then
\begin{subequations}
\begin{align}
\mu\,\log\mathsf{B}\!  &  \left(  _{\!_{\!_{\,}}}\varphi(1/2\sigma
^{2})\right)  \ \leq\ \liminf_{n\uparrow\infty}\,\,\varepsilon_{n}^{-2\beta
}m^{-\beta n}\log\mathbf{P}(R_{n}\geq\varepsilon_{n})\label{LB'}\\
\  &  \leq\ \limsup_{n\uparrow\infty}\,\varepsilon_{n}^{-2\beta}m^{-\beta
n}\log\mathbf{P}(R_{n}\geq\varepsilon_{n})\ \leq\ \mu^{-1}\log\mathsf{B}%
\!\left(  _{\!_{\!_{\,}}}\varphi(1/2\sigma^{2})\right)  \!. \label{UB'}%
\end{align}
If, additionally, $\varepsilon_{n}=m^{-\lambda_{n}/2}$ for integers
$\,\lambda_{n}\rightarrow\infty$ with $\lambda_{n}=o(n)$ as $n\uparrow\infty$,
then
\end{subequations}
\begin{equation}
\lim_{n\uparrow\infty}\varepsilon_{n}^{-2\beta}m^{-\beta n}\log\mathbf{P}%
(R_{n}\geq\varepsilon_{n})\ =\ \log\mathsf{B}\!\left(  _{\!_{\!_{\,}}}%
\varphi(1/2\sigma^{2})\right)  \!. \label{BottAsymp}%
\end{equation}
{}
\end{theorem}

According to this theorem, the main contribution to $\mathbf{P}(R_{n}%
\geq\varepsilon_{n})$ comes from lower deviations of $\,Z_{n}$ and large
deviations of $\,S_{n\,}$. In order to explain this heuristically, we note
that by Proposition~\ref{P.Boettcher} there exist (positive and finite)
constants $c_{1}\geq c_{2}$ such that
\begin{equation}
\exp\bigl[-c_{1}\,(k/m^{n})^{-\beta/(1-\beta)}\bigr]\,\leq\,m^{n}%
\,\mathbf{P}(Z_{n}=k)\,\leq\,\exp\bigl[-c_{2}\,(k/m^{n})^{-\beta/(1-\beta
)}\bigr]. \label{BBounds}%
\end{equation}
On the other hand (for details see the proof of Theorem~\ref{Bottcher} in
Subsection~\ref{SS.Bottcher} below),
\begin{equation}
\exp[-c_{3}\,\varepsilon_{n}^{2}k]\ \leq\ \mathbf{P}(S_{k}\geq\varepsilon
_{n}k)\ \leq\ \exp[-c_{4}\,\varepsilon_{n}^{2}k]
\end{equation}
for some $\,c_{3}\geq c_{4\,}.$\thinspace\ Then, roughly speaking,
\begin{equation}
\mathbf{P}(R_{n}\geq\varepsilon_{n})\ \sim\ m^{-n}\sum_{k=\mu^{n}}^{\infty
}\exp\bigl[-a\,(k/m^{n})^{-\beta/(1-\beta)}-b\,\varepsilon_{n}^{2}k\bigr]
\end{equation}
with $\,a,b>0.$\thinspace\ Obviously, the value of this sum is determined, in
a sense, by the maximal summand. It can now easily be seen, that the function
\begin{equation}
g(u)\,:=\,a\,(u/m^{n})^{-\beta/(1-\beta)}+b\,\varepsilon_{n}^{2}u,\qquad u>0,
\end{equation}
achieves its minimum at $\,u_{\ast}:=c\,\varepsilon_{n}^{-2(1-\beta)}%
m^{n\beta}$\thinspace\ [with $\,c$\thinspace\ we always denote a (positive,
finite) constant which might change its value from place to place], and,
consequently,
\begin{equation}
g(u_{\ast})\,=\,c\,\varepsilon_{n}^{2\beta}m^{n\beta}.
\end{equation}
This is in line with the normalizing sequence in Theorem~\ref{Bottcher}
(except a constant factor).

If we put formally $\alpha=\infty$ in the conditions in Theorem~\ref{T_LDev}%
(b), then (\ref{LargeDev}) should hold under the condition $\varepsilon
_{n}m^{n/2}\rightarrow\infty$, since $\varkappa\rightarrow1/2$ as
$\alpha\uparrow\infty$. But we prove it only under a slightly stronger
condition on $\,\varepsilon_{n}:$

\begin{theorem}
[\textbf{B\"{o}ttcher under heavier tails concerning }$X_{1}$]\label{T_LDev1}%
\hfill Suppose the \newline B\"{o}tt\-cher case and that $\,X_{1}$%
\thinspace\ has a tail of index $\,\theta>2$. If $\,\varepsilon_{n}%
m^{n/2}n^{-1/2\beta}\rightarrow\infty,$\thinspace\ then\/
\emph{(\ref{LargeDev})} is true.
\end{theorem}

There is the same \textquotedblleft philosophy\textquotedblright\ behind
Theorem~\ref{T_LDev1}\ as it is behind Theorem~\ref{T_LDev}(b). The main
influence of normal deviations explains also the independence of
(\ref{LargeDev}) of the parameters $\alpha$ and $\beta$. Note also that in the
special case $\,\varepsilon_{n}\equiv\varepsilon,$\thinspace
\ Theorem~\ref{T_LDev}(b) was proved in \cite{NeyVidyashankar2003}.

\begin{remark}
[\textbf{Possible generalizations}]\label{R_Conditions} Many conditions in our
results are too restrictive, but allow us to make proofs slightly shorter and
clearer. Here we mention some (almost evident) generalizations of our theorems.

\begin{itemize}
\item[(a)] It is possible to prove versions of Theorems~\ref{T_DDev} and
\ref{T_LDev} for $X_{1}$ from the domain of attraction of a stable law of any index.

\item[(b)] Theorems~\ref{T_LDev} and \ref{T_LDev1} can be generalized to the
case $\,\mathbf{P}(X_{1}\geq x)=L(x)x^{-\theta}$ with some $L,$ slowly varying
at infinity.

\item[(c)] We conjecture that condition $\,\mathbf{E}Z_{1}\log Z_{1}<\infty
$\thinspace\ can be dropped in all of our theorems. In fact, we need it only
for inequality (\ref{1.1}) below, taken from Theorem~II.4.2 of Athreya and Ney
\cite{AthreyaNey1972}. But it should be possible to prove this bound for all
supercritical Galton-Watson processes.

\item[(d)] In \cite{NeyVidyashankar2004}, $\,\mathbf{P}\bigl(Z_{n}%
\geq\varepsilon_{n\,};\,Z_{n}\geq v_{n}\bigr)$ is considered with
$\,v_{n}\rightarrow\infty$\thinspace\ and $\,\varepsilon_{n}\equiv
\varepsilon.$\thinspace\ Our methods allow to deal with the case
$\,v_{n}=o(m^{n})$\thinspace\ and $\,\varepsilon_{n}\rightarrow0$%
.\hfill$\Diamond$
\end{itemize}
\end{remark}

\begin{remark}
[\textbf{On critical Galton-Watson processes}]For the moment, suppose that the
Galton-Watson process $Z$ is critical, that is, $m=1$. Furthermore, assume
that $\,B:=\mathbf{E}Z_{1}^{2}-1\in(0,\infty).$\thinspace\ Then, analogously
to (\ref{Fluct1'}),
\begin{equation}
\lim_{n\uparrow\infty}\mathbf{P}\bigl(n^{1/2}R_{n}<x\,\big|\,Z_{n}%
>0\bigr)\ =\ \frac{2}{B}\int_{0}^{\infty}\Phi\bigl(\tfrac{xu^{1/2}}{\sigma
}\bigr)\,\mathrm{e}^{-2u/B}\,\mathrm{d}u. \label{CritFluct}%
\end{equation}
For the proof of this convergence in the two special cases of $X_{1}$ as
mentioned in Subsection~\ref{SS.motive}, see \cite{NagaevAV1967} and
\cite{Pakes1975}, respectively. {F}rom (\ref{CritFluct}) we find that for
critical processes the domain of large deviations is defined by the relation
$\varepsilon_{n}^{2}n\rightarrow\infty$ as $n\uparrow\infty$. The special case
$\varepsilon_{n}\equiv\varepsilon$ was treated by Athreya and Vidyashankar
\cite{AthreyaVidyashankar1997}. If now $\,\varepsilon_{n}\rightarrow0$ and
$\varepsilon_{n}^{2}n\rightarrow\infty$, then
\begin{equation}
\lim_{n\uparrow\infty}\varepsilon_{n}^{2}n\,\mathbf{P}\bigl(R_{n}%
\geq\varepsilon_{n}\,\big|\,Z_{n}>0\bigr)\ =\ \frac{\sigma^{2}}{B}\,.
\label{CDev2}%
\end{equation}
Actually, (\ref{CDev2}) is similar to the statement of Theorem~\ref{T_DDev} in
the case $\alpha=1$ and if $\,m^{n}$ replaced by the order $n$ of
$\,\mathbf{E}\{Z_{n}\,|\,Z_{n}>0\}$. Also, the proof of (\ref{CDev2}) is close
to the proof of Theorem~\ref{T_DDev} in the case $\alpha=1.$ There are only
two differences. First, instead of (\ref{1.1}) below, we have to use
$\,\mathbf{P}\bigl(Z_{n}=k\,\big|\,Z_{n}>0\bigr)\leq c$\thinspace$n^{-1}%
,$\thinspace\ which is derived in S.V.\ Nagaev and Wachtel
\cite{NagaevWachtel2005}. Second, we have to use the local limit theorem for
critical Galton-Watson processes instead of Proposition~\ref{P.Schroeder}. For
the proof of this local limit theorem under a second moment assumption, see
\cite{NagaevWachtel2005}.\hfill$\Diamond$
\end{remark}

\section{Auxiliary results}

In this section we prepare for the proofs of our main results.

\subsection{Separate considerations}

As a first step, we state two bounds for local probabilities of our
supercritical Galton-Watson process $\,Z$\thinspace\ (satisfying
$\,\mathbf{E}Z_{1}\log Z_{1}<\infty).$

\begin{lemma}
[\textbf{Local probabilities of} \thinspace$Z$]\label{L1}We have%
\begin{equation}
\mathbf{P}\bigl(Z_{n}=k\,\big|\,Z_{0}=\ell\bigr)\ \leq\ c\ \frac{\ell}%
{k}\,,\qquad k,\ell,n\geq1. \label{1.1}%
\end{equation}
Moreover, in the Schr\"{o}der case,%
\begin{equation}
\mathbf{P}\bigl(Z_{n}=k\,\big|\,Z_{0}=1\bigr)\ \leq\ c\ \frac{k^{\alpha-1}%
}{m^{\alpha n}}\,,\qquad k,n\geq1. \label{1.2}%
\end{equation}
{}
\end{lemma}

\begin{proof}
For aperiodic $(d=1)$ offspring laws the proof of inequality (\ref{1.1}) is
given in \cite[Theorem~II.4.2]{AthreyaNey1972}. The proof in the remaining
case $d>1$ can be carried out similarly.

In proving (\ref{1.2}) it is sufficient to assume that $k\leq m^{n}$,
otherwise (\ref{1.2}) follows from (\ref{1.1}). Under the present condition
$\,\mathbf{E}Z_{1}\log Z_{1}<\infty,$\thinspace\ formula (151) from
\cite{FleischmannWachtel2006.LowerDev.AIHP} with $\,N=\ell_{0}:=1+[1/\alpha
]$\thinspace\ and $\,j=n-a_{k}$\thinspace\ where $\,a_{k}:=\min\{j\geq
1:\,m^{j}\geq k\}$\thinspace\ reads as
\begin{equation}
\sum_{\ell=\ell_{0}}^{\infty}\mathbf{P}(Z_{n-a_{k}}=\ell)\,\mathbf{P}%
\bigl(Z_{a_{k}}=k\,\big|\,Z_{0}=\ell\bigr)\ \leq\ \frac{c}{m^{a_{k}}%
}\,f_{n-a_{k}}(\mathrm{e}^{-\delta}).
\end{equation}
It follows from (\ref{SConv'}) that the right hand side is bounded by
$c$\thinspace$m^{-a_{k}}\gamma^{n-a_{k}}$. Since%
\begin{equation}
k\,\leq\,m^{a_{k}}\,\leq\,mk\quad\text{and}\quad\gamma=m^{-\alpha},
\label{again}%
\end{equation}
we get the bound
\begin{equation}
\sum_{\ell=\ell_{0}}^{\infty}\mathbf{P}(Z_{n-a_{k}}=\ell)\,\mathbf{P}%
\bigl(Z_{a_{k}}=k\,\big|\,Z_{0}=\ell\bigr)\ \leq\ c\ \frac{k^{\alpha-1}%
}{m^{\alpha n}}\,. \label{1.3}%
\end{equation}
If $\,\ell_{0}=1$, then the proof of (\ref{1.2}) is complete, since the left
hand side in (\ref{1.3}) equals $\mathbf{P}(Z_{n}=k)$. Assume now that
$\ell_{0}\geq2$. From (\ref{1.1}) it follows that
\begin{equation}
\sum_{\ell=1}^{\ell_{0}-1}\mathbf{P}(Z_{n-a_{k}}=\ell)\,\mathbf{P}%
\bigl(Z_{a_{k}}=k\,\big|\,Z_{0}=\ell\bigr)\ \leq\ c\ \frac{\ell_{0}}{k}%
\,\sum_{\ell=1}^{\ell_{0}-1}\mathbf{P}(Z_{n-a_{k}}=\ell).
\end{equation}
By (\ref{SConv}), $\,\lim_{n\uparrow\infty}\gamma^{-n}\,\mathbf{P}(Z_{n}%
=\ell)=\nu_{\ell}<\infty,$\thinspace\ for every fixed $\ell.$ Hence,
\begin{equation}
\sum_{\ell=1}^{\ell_{0}-1}\mathbf{P}(Z_{n-a_{k}}=\ell)\ \leq\ c\,\gamma
^{n-a_{k}}%
\end{equation}
for all $n\geq1$. Using again (\ref{again}), we get
\begin{equation}
\sum_{\ell=1}^{\ell_{0}-1}\mathbf{P}(Z_{n-a_{k}}=\ell)\,\mathbf{P}%
\bigl(Z_{a_{k}}=k\,\big|\,Z_{0}=\ell\bigr)\ \leq\ c\ \frac{k^{\alpha-1}%
}{m^{\alpha n}}\,.
\end{equation}
This completes the proof.
\end{proof}

The following lemma contains two versions of the so-called Fuk-Nagaev
inequality for tail probabilities of sums of i.i.d.\ variables. Recall that we
assumed that $X_{1}$ is centered and has positive finite variance.

\begin{lemma}
[\textbf{Fuk-Nagaev inequality}]\label{L2}For $\,k\geq1,\,\ \varepsilon
_{n}>0,$\thinspace\ $n\geq1,\,\ r>1,\,\ $and $\,t\geq2$,
\begin{equation}
\mathbf{P}(S_{k}\geq\varepsilon_{n}k)\ \leq\ k\,\mathbf{P}\bigl(X_{1}\geq
r^{-1}\varepsilon_{n}k\bigr)+(\mathrm{e}\,r\sigma^{2})^{r}\,\varepsilon
_{n}^{-2r}\,k^{-r}, \label{1.4}%
\end{equation}
and
\begin{align}
\mathbf{P}(S_{k}\geq\varepsilon_{n}k)\ \leq &  \ k\,\mathbf{P}\bigl(X_{1}\geq
r^{-1}\varepsilon_{n}k\bigr)+\exp\!\Big[-\frac{2}{(t+2)^{2}\,\mathrm{e}%
^{t}\,\sigma^{2}}\ \varepsilon_{n}^{2}k\Big ]\label{1.4'}\\
&  +\Bigl(\frac{(t+2)\,r^{t-1}\,\mathbf{E}\!\left\{  _{\!_{\!_{\,}}}X_{1}%
^{t};\ 0\leq X_{1}\leq\varepsilon_{n}k\right\}  }{t\,\varepsilon_{n}%
^{t}\,k^{t-1}}\Bigr)^{tr/(t+2)}.\nonumber
\end{align}
{}
\end{lemma}

\begin{proof}
By (1.56) and (1.23) in S.V.\ Nagaev \cite{NagaevSV1979}, for all $\,u,v>0,$
\begin{equation}
\mathbf{P}(S_{k}\geq u)\ \leq\ k\,\mathbf{P}(X_{1}\geq v)+\mathrm{e}%
^{u/v}\Bigl(\frac{\sigma^{2}k}{uv}\Bigr)^{\!u/v}%
\end{equation}
and
\begin{align}
\mathbf{P}(S_{k}\geq u)\ \leq &  \ k\,\mathbf{P}(X_{1}\geq v)+\exp
\!\Big[-\frac{2u^{2}}{(t+2)^{2}\,\mathrm{e}^{t}\,\sigma^{2}}\Big ]\\
&  +\Bigl(\frac{(t+2)\,k\,\mathbf{E}\!\left\{  _{\!_{\!_{\,}}}X_{1}%
^{t};\ 0\leq X_{1}\leq v\right\}  }{t\,u\,v^{t-1}}\Bigr)^{\!tu/(t+2)v}%
.\nonumber
\end{align}
Putting here $\,u=\varepsilon_{n}k$\thinspace\ and $\,v=u/r,$\thinspace\ we
get (\ref{1.4}) and (\ref{1.4'}), finishing the proof.
\end{proof}

\begin{remark}
[\textbf{On the case $\varepsilon_{n}\equiv\varepsilon$}]Here we prove a
one-sided version of~(\ref{Athreya1})~con\-cerning our general $\,R_{n\,}%
,$\thinspace\ assuming the Schr\"{o}der case and that $\,\mathbf{E}(X_{1}%
^{+})^{1+\alpha}<\infty.$\thinspace\ Take any $\varepsilon>0$ and set
$\,g_{n}(k):=m^{\alpha n}$\thinspace$\mathbf{P}(Z_{n}=k)\,\mathbf{P}(S_{k}%
\geq\varepsilon k)$. From estimate (\ref{1.2}) we get, for all $\,n,k\geq
1,$\thinspace\ the inequality $\,g_{n}(k)\leq c$\thinspace$\tilde{g}%
(k),$\thinspace\ where $\,\tilde{g}(k):=k^{\alpha-1}$\thinspace$\mathbf{P}%
(S_{k}\geq\varepsilon k).$\thinspace\ Next we show that $\,\tilde{g}%
(k)$\thinspace\ is summable in $k$. Letting $\,\varepsilon_{n}=\varepsilon
$\thinspace\ and $\,r=\alpha+1$\thinspace\ in (\ref{1.4}), we see that for all
$\,k\geq1$,
\begin{equation}
\tilde{g}(k)\ \leq\ k^{\alpha}\,\mathbf{P}\bigl(X_{1}\geq\varepsilon
k/(1+\alpha)\bigr)+c\,\varepsilon^{-2-2\alpha}k^{-2}.
\end{equation}
But the summability of $\,k^{\alpha}$\thinspace$\mathbf{P}(X_{1}\geq ck)$ with
some (hence all) positive $\,c$\thinspace\ is equivalent to the finiteness of
$\,\mathbf{E}(X_{1}^{+})^{1+\alpha}$, and we get the claimed summability of
$\,\tilde{g}(k)$.

On the other hand, it follows from (\ref{SConv}) that for every fixed $k$,
\begin{equation}
\lim_{n\uparrow\infty}g_{n}(k)\ =\ \nu_{k}\,\mathbf{P}(S_{k}\geq\varepsilon
k).
\end{equation}
Therefore, by dominated convergence,
\begin{equation}
\lim_{n\uparrow\infty}\sum_{k=1}^{\infty}g_{n}(k)\ =\ \sum_{k=1}^{\infty}%
\nu_{k}\,\mathbf{P}(S_{k}\geq\varepsilon k).
\end{equation}
Recalling the definition of $g_{n}(k)$ and using (\ref{TotalProb}), we obtain%
\begin{equation}
\lim_{n\uparrow\infty}\,m^{\alpha n}\,\mathbf{P}(R_{n}\geq\varepsilon
_{n})\ =\ \sum_{k=1}^{\infty}\nu_{k}\,\mathbf{P}(S_{k}\geq\varepsilon k),
\end{equation}
yielding the wanted one-sided version.\hfill$\Diamond$
\end{remark}

\subsection{Interplay between the two competing forces}

In the next four lemmas we prove bounds for different parts of the sum at the
right hand side of (\ref{TotalProb}).

\begin{lemma}
[\textbf{A tail estimate}]\label{L4}Assume $\,X_{1}$\thinspace\ has a tail of
index $\,\theta>2$.\thinspace\ Then%
\begin{align}
\sum_{k\,\geq m^{n}}  &  \mathbf{P}(Z_{n}=k)\,\mathbf{P}(S_{k}\geq
\varepsilon_{n}k)\\
\  &  \leq\ c\,\Bigl(\varepsilon_{n}^{-\theta}\,m^{-(\theta-1)n}%
+(\varepsilon_{n}^{2}m^{n})^{-1}\exp\bigl[-c\,\varepsilon_{n}^{2}%
\,m^{n}\bigr]\Bigr),\qquad\varepsilon_{n}>0,\ \,n\geq1.\nonumber
\end{align}

\end{lemma}

\begin{proof}
Letting $t=\theta+1$ and $r=(t+2)/t$\thinspace\ in (\ref{1.4'}), and using
that $\,X_{1}$\thinspace\ has a tail of index $\,\theta>2,$\thinspace\ we get
the bound
\begin{equation}
\mathbf{P}(S_{k}\geq\varepsilon_{n}k)\ \leq\ c\,\Bigl(\varepsilon_{n}%
^{-\theta}\,k^{-(\theta-1)}+\frac{\mathbf{E}\bigl\{X_{1}^{\theta+1};\,X_{1}%
\in\lbrack0,\varepsilon_{n}k]\bigr\}}{\varepsilon_{n}^{\theta+1}\,k^{\theta}%
}\Bigr)+\exp[-c\,\varepsilon_{n}^{2}k]. \label{57}%
\end{equation}
Clearly, under (\ref{TailAssumption}),
\begin{equation}
\mathbf{E}\bigl\{X_{1}^{\theta+1};\ X_{1}\in\lbrack0,x]\bigr\}\,\sim
\,a\,\theta x\quad\text{as }\,x\uparrow\infty. \label{mom.tail}%
\end{equation}
Thus,
\begin{equation}
\mathbf{E}\bigl\{X_{1}^{\theta+1};\ X_{1}\in\lbrack0,x]\bigr\}\ \leq
\ c\,x,\qquad x\geq1.
\end{equation}
On the other hand, if $\,x\leq1,$\thinspace\
\begin{equation}
\mathbf{E}\bigl\{X_{1}^{\theta+1};\ X_{1}\in\lbrack0,x]\bigr\}\ \leq
\ x^{\theta+1}\,\mathbf{P}\!\left(  X_{1}\in\lbrack0,x]\right)  \ \leq\ x.
\end{equation}
Therefore,%
\begin{equation}
\mathbf{E}\bigl\{X_{1}^{\theta+1};\ X_{1}\in\lbrack0,x]\bigr\}\ \leq
\ c\,x,\qquad x\geq0.
\end{equation}
Applying this to the expectation in (\ref{57}), we get
\begin{equation}
\mathbf{P}(S_{k}\geq\varepsilon_{n}k)\ \leq\ c\,\varepsilon_{n}^{-\theta
}\,k^{-(\theta-1)}+\exp[-c\,\varepsilon_{n}^{2}\,k].
\end{equation}
Moreover, combining this bound with (\ref{1.1}) gives%
\begin{equation}
\sum_{k\,\geq m^{n}}\mathbf{P}(Z_{n}=k)\,\mathbf{P}(S_{k}\geq\varepsilon
_{n}k)\ \leq\ c\,\varepsilon_{n}^{-\theta}\sum_{k\,\geq m^{n}}k^{-\theta}%
+\sum_{k\,\geq m^{n}}k^{-1}\exp[-c\,\varepsilon_{n}^{2}k]. \label{1.9}%
\end{equation}
Obviously,
\begin{equation}
\sum_{k\,\geq m^{n}}k^{-\theta}\ \leq\ c\,m^{-(\theta-1)n}. \label{1.10}%
\end{equation}
On the other hand,
\begin{gather}
\sum_{k\,\geq m^{n}}k^{-1}\exp[-c\,\varepsilon_{n}^{2}k]\ \leq\ m^{-n}%
\sum_{k\,\geq m^{n}}\exp[-c\,\varepsilon_{n}^{2}k]\nonumber\\
\ \leq\ c\,(\varepsilon_{n}^{2}m^{n})^{-1}\exp\bigl[-c\,\varepsilon_{n\,}%
^{2}m^{n}\bigr]. \label{1.11}%
\end{gather}
Substituting (\ref{1.10}) and (\ref{1.11}) into (\ref{1.9}) finishes the proof.
\end{proof}

\begin{lemma}
[\textbf{Another tail estimate}]\label{L5}Suppose the Schr\"{o}der case, let
$\,X_{1}$\thinspace\ satisfy the moment condition\/ \emph{(\ref{MomCond}),}
and let $\,\varepsilon_{n}\rightarrow0.$\thinspace\ Then
\begin{equation}
\limsup_{n\uparrow\infty}\,\varepsilon_{n}^{2\alpha}\,m^{\alpha n}%
\sum_{k\,\geq A/\varepsilon_{n}^{2}}\mathbf{P}(Z_{n}=k)\,\mathbf{P}(S_{k}%
\geq\varepsilon_{n}k)\ \leq\ \frac{c}{A}\,,\qquad A\geq1.
\end{equation}
{}
\end{lemma}

\begin{proof}
Combining (\ref{1.2}) and (\ref{1.4}) with $\,r=\alpha+1$\thinspace\ gives
\begin{align}
&  m^{\alpha n}\sum_{k\,\geq A/\varepsilon_{n}^{2}}\mathbf{P}(Z_{n}%
=k)\,\mathbf{P}(S_{k}\geq\varepsilon_{n}k)\nonumber\\
\  &  \leq\ c\,\Bigl(\sum_{k\,\geq A/\varepsilon_{n}^{2}}k^{\alpha
}\,\mathbf{P}\bigl(X_{1}\geq(\alpha+1)^{-1}\varepsilon_{n}%
k\bigr)\,+\,\varepsilon_{n}^{-2(\alpha+1)}\sum_{k\,\geq A/\varepsilon_{n}^{2}%
}k^{-2}\Bigr). \label{1.12}%
\end{align}
Note that
\begin{equation}
\varepsilon_{n}^{-2(\alpha+1)}\Big(\sum_{k\,\geq A/\varepsilon_{n}^{2}}%
k^{-2}\Big)\ \leq\ \frac{c}{A}\,\varepsilon_{n}^{-2\alpha},\qquad
n>0,\ \,\varepsilon_{n}>0,\,\ A\geq1. \label{1.13}%
\end{equation}
On the other hand, to bound the first sum at the right hand side in
(\ref{1.12}), note first that%
\[
\int_{k-1}^{k}u^{\alpha}\,\mathbf{P}(X_{1}\geq(\alpha+1)^{-1}\varepsilon
_{n}u)\,\mathrm{d}u\ \geq\ (k-1)^{\alpha}\,\mathbf{P}\!\left(  _{\!_{\!_{\,}}%
}X_{1}\geq(\alpha+1)^{-1}\varepsilon_{n}k\right)  \!,\qquad k\geq1.
\]
This inequality can be continued by using $\,k-1\geq k/2$\thinspace\ for
$\,k\geq2.$\thinspace\ Summing up gives for $\,\varepsilon_{n}^{2}\leq1/2$,
\begin{gather}
\sum_{k\,\geq A/\varepsilon_{n}^{2}}k^{\alpha}\,\mathbf{P}\bigl(X_{1}%
\geq(\alpha+1)^{-1}\varepsilon_{n}k\bigr)\leq\ c\int_{A/\varepsilon_{n}^{2}%
-1}^{\infty}u^{\alpha}\,\mathbf{P}\!\left(  _{\!_{\!_{\,}}}X_{1}\geq
(\alpha+1)^{-1}\varepsilon_{n}u\right)  \mathrm{d}u\nonumber\\
\ \leq\ c\,\varepsilon_{n}^{-\alpha-1}\int_{(A-\varepsilon_{n}^{2}%
)/(\alpha+1)\varepsilon_{n}}^{\infty}v^{\alpha}\,\mathbf{P}(X_{1}\geq
v)\,\mathrm{d}v. \label{1.14}%
\end{gather}

Recall that we assumed the moment condition (\ref{MomCond}) and that
$\,\varepsilon_{n}\rightarrow0.$\thinspace\ Then the integral in (\ref{1.14})
converges to zero as $\,n\uparrow\infty,$\thinspace\ uniformly in $\,A\geq
1.$\thinspace\ In particular, under $\,\alpha\geq1,$\thinspace\ (\ref{1.14})
is of order $\,o(\varepsilon_{n}^{-2\alpha}),$\thinspace\ uniformly in
$\,A\geq1.$\thinspace\ On the other hand, if $\,\alpha<1$\thinspace\ and since
$\,\mathbf{E}X_{1}^{2}<\infty$,%
\begin{align}
&  \int_{(A-\varepsilon_{n}^{2})/(\alpha+1)\varepsilon_{n}}^{\infty}v^{\alpha
}\,\mathbf{P}(X_{1}\geq v)\,\mathrm{d}v\ \\
&  \leq\ c\ \frac{\varepsilon_{n}^{1-\alpha}}{(A-\varepsilon_{n}%
^{2})^{1-\alpha}}\int_{(A-\varepsilon_{n}^{2})/(\alpha+1)\varepsilon_{n}%
}^{\infty}v\,\mathbf{P}(X_{1}\geq v)\,\mathrm{d}v\ =\ o(\varepsilon
_{n}^{1-\alpha})\ =\ o(\varepsilon_{n}^{-2\alpha})\nonumber
\end{align}
as $\,n\uparrow\infty,$\thinspace\ uniformly in $\,A\geq1.$\thinspace\ Thus,
for each $\,\alpha<\infty$\thinspace\ we have
\begin{equation}
\sup_{A\geq1}\ \sum_{k\,\geq A/\varepsilon_{n}^{2}}k^{\alpha}\,\mathbf{P}%
\bigl(X_{1}\geq(\alpha+1)^{-1}\varepsilon_{n}k\bigr)\ =\ o(\varepsilon
_{n}^{-2\alpha})\quad\text{as}\,\ n\uparrow\infty. \label{1.15}%
\end{equation}
In particular,%
\begin{equation}
\limsup_{n\uparrow\infty}\,\varepsilon_{n}^{2\alpha}\sum_{k\,\geq
A/\varepsilon_{n}^{2}}k^{\alpha}\,\mathbf{P}\bigl(X_{1}\geq(\alpha
+1)^{-1}\varepsilon_{n}k\bigr)\ \leq\ \frac{c}{A}\,,\qquad A\geq1.
\label{1.15'}%
\end{equation}
Combining (\ref{1.12}), (\ref{1.13}), and (\ref{1.15'}) gives the claim in the lemma.
\end{proof}

\begin{lemma}
[\textbf{Initial part}]\label{L6'}In the Schr\"{o}der case,%
\begin{equation}
\sum_{1\,\leq\,k\,\leq\,\delta/\varepsilon_{n}^{2}}\mathbf{P}(Z_{n}%
=k)\,\mathbf{P}(S_{k}\geq\varepsilon_{n}k)\ \leq\ c\,\delta^{\alpha
}\varepsilon_{n}^{-2\alpha}\,m^{-\alpha n}, \label{1.16}%
\end{equation}
$\delta>0,\ \,\varepsilon_{n}>0,\ \,n\geq1.$
\end{lemma}

\begin{proof}
It follows from (\ref{1.2}) that
\begin{gather}
\sum_{1\,\leq\,k\,\leq\,\delta/\varepsilon_{n}^{2}}\mathbf{P}(Z_{n}%
=k)\,\mathbf{P}(S_{k}\geq\varepsilon_{n}k)\ \leq\ \sum_{1\,\leq\,k\,\leq
\,\delta/\varepsilon_{n}^{2}}\mathbf{P}(Z_{n}=k)\ \nonumber\\
\leq\ \frac{c}{m^{\alpha n}}\,\sum_{1\,\leq\,k\,\leq\,\delta/\varepsilon
_{n}^{2}}k^{\alpha-1}\ \leq\ c\,\delta^{\alpha}\varepsilon_{n}^{-2\alpha
}\,m^{-\alpha n},
\end{gather}
finishing the proof.
\end{proof}

\begin{lemma}
[\textbf{A central part and another initial part estimate}]\label{L6}\hfill
Suppose $\,\newline1<\alpha<\infty$\thinspace\ and that $\,X_{1}$%
\thinspace\ has a tail of index $\,\theta\in(2,1+\alpha)$.\thinspace\ Then%
\begin{align}
\sum_{A/\varepsilon_{n}^{2}\,\leq\,k\,\leq\,\delta m^{n}}  &  \mathbf{P}%
(Z_{n}=k)\,\mathbf{P}(S_{k}\geq\varepsilon_{n}k)\ \label{6.1}\\[1pt]
&  \leq\ c\,\bigl(\delta^{1+\alpha-\theta}\varepsilon_{n}^{-\theta}%
m^{-(\theta-1)n}+A^{-1}\,\varepsilon_{n}^{-2\alpha}\,m^{-\alpha n}%
\bigr),\nonumber
\end{align}
$A\geq1,\ \,\delta>0,\,\ \varepsilon_{n}>0,\,\ n\geq1,$\thinspace\ and{}%
\begin{align}
&  \sum_{1\,\leq\,k\,\leq\,\delta m^{n}}\mathbf{P}(Z_{n}=k)\,\mathbf{P}%
(S_{k}\geq\varepsilon_{n}k)\ \label{6.2}\\
&  \leq\ c\,\bigl(\delta^{1+\alpha-\theta}\varepsilon_{n}^{-\theta}%
m^{-(\theta-1)n}+\varepsilon_{n}^{-2\alpha}\,m^{-\alpha n}\bigr),\qquad
\delta>0,\,\ \varepsilon_{n}>0,\,\ n\geq1.\nonumber
\end{align}

\end{lemma}

\begin{proof}
Combining (\ref{1.2}) and (\ref{1.4}) with $\,r=\alpha+1$ gives
\begin{align}
&  \sum_{A/\varepsilon_{n}^{2}\,\leq\,k\,\leq\,\delta m^{n}}\mathbf{P}%
(Z_{n}=k)\,\mathbf{P}(S_{k}\geq\varepsilon_{n}k)\ \label{74}\\
&  \leq\,c\,m^{-\alpha n}\Bigg(\sum_{A/\varepsilon_{n}^{2}\,\leq
\,k\,\leq\,\delta m^{n}}\!k^{\alpha}\,\mathbf{P}\bigl(X_{1}\geq(\alpha
+1)^{-1}\varepsilon_{n}k\bigr)\!+\varepsilon_{n}^{-2(\alpha+1)}\!\!\!\sum
_{A/\varepsilon_{n}^{2}\,\leq\,k\,\leq\,\delta m^{n}}\!k^{-2}\!\Bigg)\hspace
{-1pt}.\nonumber
\end{align}
From (\ref{1.13}),
\begin{equation}
\varepsilon_{n}^{-2(\alpha+1)}\sum_{A/\varepsilon_{n}^{2}\,\leq\,k\,\leq
\,\delta m^{n}}k^{-2}\ \leq\ \frac{c}{A}\,\varepsilon_{n}^{-2\alpha}.
\label{1.17}%
\end{equation}
On the other hand, since $\,X_{1}$\thinspace\ has a tail of index $\,\theta
\in(2,1+\alpha),$
\begin{gather}
\sum_{A/\varepsilon_{n}^{2}\,\leq\,k\,\leq\,\delta m^{n}}k^{\alpha
}\,\mathbf{P}\bigl(X_{1}\geq(\alpha+1)^{-1}\varepsilon_{n}k\bigr)\ \leq
\ c\,\varepsilon_{n}^{-\theta}\sum_{1\,\leq\,k\,\leq\,\delta m^{n}}%
k^{\alpha-\theta}\nonumber\\
\ \leq\ c\,\varepsilon_{n}^{-\theta}\delta^{1+\alpha-\theta}m^{(1+\alpha
-\theta)n}. \label{1.18}%
\end{gather}
Combine (\ref{74})--(\ref{1.18}) to get (\ref{6.1}).

Putting $\,A=1$ in (\ref{6.1}) and $\delta=1$ in (\ref{1.16}), we obtain
(\ref{6.2}), finishing the proof.
\end{proof}

Recall that $(\mu,d)$ refers to the type of the offspring law, $\,\alpha
\in(0,\infty)$\thinspace\ to the Schr\"{o}der constant, and that $\,X_{1}%
$\thinspace\ is assumed to have a finite variance. For $\,0<\delta
<1<A<\infty,$\thinspace\ consider
\begin{equation}
\Sigma_{n}(\delta,A)\ :=\ \sum_{\delta/\varepsilon_{n}^{2}\,\leq
\,k\,\leq\,A/\varepsilon_{n}^{2}}\mathbf{P}(Z_{n}=k)\,\mathbf{P}(S_{k}%
\geq\varepsilon_{n}k).
\end{equation}

\begin{lemma}
[\textbf{Another central part estimate}]\label{L7}Suppose the Schr\"{o}der
case, that $\,\varepsilon_{n}\rightarrow0,$\thinspace\ and that $\,\varepsilon
_{n}^{2}m^{n}\rightarrow\infty$.\thinspace\ Then for all $\,0<\delta
<1<A<\infty,$
\begin{equation}%
\begin{array}
[c]{l}%
\displaystyle
V_{\ast}\int_{\delta}^{A}u^{\alpha-1}\ \overline{\Phi}(\sqrt{u}/\sigma
)\,\mathrm{d}u\ \leq\ \liminf_{n\uparrow\infty}\,\varepsilon_{n}^{2\alpha
}\,m^{\alpha n}\,\Sigma_{n}(\delta,A)\\%
\displaystyle
\leq\ \limsup_{n\uparrow\infty}\,\varepsilon_{n}^{2\alpha}\,m^{\alpha
n}\,\Sigma_{n}(\delta,A)\ \leq\ V^{\ast}\!\!\int_{\delta}^{A}u^{\alpha
-1}\ \overline{\Phi}(\sqrt{u}/\sigma)\,\mathrm{d}u
\end{array}
\label{22}%
\end{equation}
with $\,V_{\ast}$\thinspace\ and $\,V^{\ast}$\thinspace\ defined in\/
\emph{(\ref{not.V.star})}, and where $\,\overline{\Phi}(x):=1-\Phi(x).$
\end{lemma}

\begin{proof}
In view of (\ref{44}) in Proposition~\ref{P.Schroeder} with $\,k_{n}%
=\delta/\varepsilon_{n\,}^{2},$\thinspace\
\begin{equation}
\Sigma_{n}(\delta,A)=(1+o(1))d\sum_{k\in H(\delta,A)}m^{-n}\,w\bigl(\tfrac
{k}{m^{n}}\bigr)\,\mathbf{P}(S_{k}\geq\varepsilon_{n}k)\quad\text{as}%
\,\ n\uparrow\infty
\end{equation}
with $H(\delta,A):=\bigl\{k\in\lbrack\delta/\varepsilon_{n}^{2}%
,\,A/\varepsilon_{n}^{2}]:\,k\equiv\mu\ (\mathrm{mod}\,d)\bigr\}.$%
\thinspace\ Clearly,
\begin{gather}
V_{\ast}(n)\sum_{k\in H(\delta,A)}\frac{k^{\alpha-1}}{m^{\alpha n}%
}\,\mathbf{P}(S_{k}\geq\varepsilon_{n}k)\ \leq\ \sum_{k\in H(\delta,A)}%
m^{-n}\,w\bigl(\tfrac{k}{m^{n}}\bigr)\,\mathbf{P}(S_{k}\geq\varepsilon
_{n}k)\label{82}\\
\ \leq\ V^{\ast}(n)\sum_{k\in H(\delta,A)}\frac{k^{\alpha-1}}{m^{\alpha n}%
}\ \mathbf{P}(S_{k}\geq\varepsilon_{n}k),\nonumber
\end{gather}
where we set
\begin{equation}
V_{\ast}(n)\ :=\ \inf_{u\leq A/\varepsilon_{n}^{2}m^{n}}u^{1-\alpha}w(u),\quad
V^{\ast}(n)\ :=\ \sup_{u\leq A/\varepsilon_{n}^{2}m^{n}}u^{1-\alpha}w(u).
\label{88}%
\end{equation}
By the central limit theorem,
\begin{equation}
\sup_{k\in H(\delta,A)}\Big|\mathbf{P}(S_{k}\geq\varepsilon_{n}k)-\overline
{\Phi}\bigl(\sqrt{\varepsilon_{n}^{2}k}/\sigma\bigr)\Big|\rightarrow
0\quad\text{as}\,\ n\uparrow\infty.
\end{equation}
Hence, as $n\uparrow\infty$,
\begin{gather}
\sum_{k\in H(\delta,A)}k^{\alpha-1}\,\mathbf{P}(S_{k}\geq\varepsilon
_{n}k)\ =\ \left(  _{\!_{\!_{\,}}}1+o(1)\right)  \sum_{k\in H(\delta
,A)}k^{\alpha-1}\ \overline{\Phi}\bigl(\sqrt{\varepsilon_{n}^{2}k}%
/\sigma\bigr)\nonumber\\
\ =\ \varepsilon_{n}^{-2\alpha}\left(  _{\!_{\!_{\,}}}1+o(1)\right)
\sum_{k\in H(\delta,A)}(\varepsilon_{n}^{2}k)^{\alpha-1}\ \overline{\Phi
}\bigl(\sqrt{\varepsilon_{n}^{2}k}/\sigma\bigr)\,\varepsilon_{n}%
^{2}\nonumber\\
\ =\ d^{-1}\varepsilon_{n}^{-2\alpha}\left(  _{\!_{\!_{\,}}}1+o(1)\right)
\int_{\delta}^{A}u^{\alpha-1}\ \overline{\Phi}(\sqrt{u}/\sigma)\,\mathrm{d}u.
\label{85}%
\end{gather}
Substituting (\ref{85}) into (\ref{82}) and noting that we have $V_{\ast
}(n)\rightarrow V_{\ast}$ and $V^{\ast}(n)\rightarrow V^{\ast}$ as
$n\uparrow\infty$ by our velocity assumption on $\,\varepsilon_{n\,}%
,$\thinspace\ we obtain (\ref{22}).
\end{proof}

Finally, we compute the limit, as $\delta\downarrow0$ and $A\uparrow\infty$,
of the integral from (\ref{22}).

\begin{lemma}
[\textbf{A moment formula for the Gaussian law}]\label{L8}For $\,0<\alpha
<\infty$,
\begin{equation}
\int_{0}^{\infty}u^{\alpha-1}\ \overline{\Phi}(\sqrt{u}/\sigma)\,\mathrm{d}%
u\ =\ \frac{2^{\alpha-1\,}\Gamma(\alpha+1/2)}{\alpha\sqrt{\pi}}\ \sigma
^{2\alpha}\ =\ \Gamma_{\alpha\,}. \label{inL8}%
\end{equation}
{}
\end{lemma}

\begin{proof}
Substituting $\,v=\sqrt{u}/\sigma,$\thinspace\ we have
\begin{align}
&  \int_{0}^{\infty}u^{\alpha-1}\ \overline{\Phi}(\sqrt{u}/\sigma
)\,\mathrm{d}u\ =\ 2\sigma^{2\alpha}\int_{0}^{\infty}v^{2\alpha-1}%
\,\overline{\Phi}(v)\,\mathrm{d}v\nonumber\\
\  &  =\ 2\sigma^{2\alpha}\int_{0}^{\infty}\!\mathrm{d}v\ v^{2\alpha-1}%
\int_{v}^{\infty}\!\mathrm{d}t\ \frac{1}{\sqrt{2\pi}}\,\mathrm{e}^{-t^{2}%
/2}\ =\ 2\ \frac{\sigma^{2\alpha}}{\sqrt{2\pi}}\int_{0}^{\infty}%
\!\mathrm{d}t\ \mathrm{e}^{-t^{2}/2}\int_{0}^{t}\!\mathrm{d}v\ v^{2\alpha
-1}\nonumber\\
\  &  =\ \frac{\sigma^{2\alpha}}{\alpha\sqrt{2\pi}}\int_{0}^{\infty}%
t^{2\alpha}\,\mathrm{e}^{-t^{2}/2}\,\mathrm{d}t.
\end{align}
Substituting now $v=t^{2}/2$, the chain of equalities can be continued with
\begin{equation}
=\ \frac{2^{\alpha-1}\sigma^{2\alpha}}{\alpha\sqrt{\pi}}\int_{0}^{\infty
}v^{\alpha-1/2}\,\mathrm{e}^{-v}\,\mathrm{d}v\ =\ \frac{2^{\alpha-1\,}%
\Gamma(\alpha+1/2)}{\alpha\sqrt{\pi}}\ \sigma^{2\alpha},
\end{equation}
which equals $\,\Gamma_{\alpha}$\thinspace\ from (\ref{notC.alpha}). The proof
is finished.
\end{proof}

\section{Proof of the main results}

\subsection{Proof of Theorem~\ref{T_DDev}}

We will use decomposition (\ref{TotalProb}). Combining Lemmas~\ref{L6'},
\ref{L7}, and \ref{L5}, and using that $\delta$ and $A$ are arbitrary, we see
that
\begin{align}
&  V_{\ast}\int_{0}^{\infty}u^{\alpha-1}\ \overline{\Phi}(\sqrt{u}%
/\sigma)\,\mathrm{d}u\ \leq\ \liminf_{n\uparrow\infty}\,\varepsilon
_{n}^{2\alpha}\,m^{\alpha n}\sum_{k=1}^{\infty}\mathbf{P}(Z_{n}=k)\,\mathbf{P}%
(S_{k}\geq\varepsilon_{n}k)\\
\  &  \leq\ \limsup_{n\uparrow\infty}\,\varepsilon_{n}^{2\alpha}\,m^{\alpha
n}\sum_{k=1}^{\infty}\mathbf{P}(Z_{n}=k)\,\mathbf{P}(S_{k}\geq\varepsilon
_{n}k)\ \leq\ V^{\ast}\!\!\int_{0}^{\infty}u^{\alpha-1}\ \overline{\Phi}%
(\sqrt{u}/\sigma)\,\mathrm{d}u.\nonumber
\end{align}
With Lemma~\ref{L8} the proof is finished.\hfill$\square$

\subsection{Proof of Theorem~\ref{T_LDev}}

Let $\,\varepsilon_{n}=o(m^{-\varkappa n}).$\thinspace\ Then under the
assumptions in the theorem,
\begin{equation}
\varepsilon_{n}^{-\theta}m^{-(\theta-1)n}+(\varepsilon_{n}^{2}m^{n}%
)^{-1}\,=\,o(\varepsilon_{n}^{-2\alpha}\,m^{-\alpha n}). \label{91}%
\end{equation}
{F}rom this relation, Lemma~\ref{L4}, and (\ref{6.1}) with $\,\delta
=1,$\thinspace\ we get%
\begin{equation}
\limsup_{n\uparrow\infty}\,\varepsilon_{n}^{2\alpha}\,m^{\alpha n}%
\sum_{k\,\geq A/\varepsilon_{n}^{2}}\mathbf{P}(Z_{n}=k)\,\mathbf{P}(S_{k}%
\geq\varepsilon_{n}k)\ \leq\ \frac{c}{A}\,,\qquad A\geq1.
\end{equation}
Theorem~\ref{T_LDev}(a) follows from this bound, (\ref{1.16}), and
(\ref{22}).\smallskip

We turn now to the proof of parts (b) and (c). It is known (see for example
Borovkov \cite{Borovkov2000}), that if $\,\mathbf{P}(X_{1}\geq x)$ is
regularly varying as $\,x\uparrow\infty$\thinspace\ with index $\theta>2$,
then for every sequence $\,a_{k}\rightarrow\infty$,
\begin{equation}
\lim_{k\uparrow\infty}\ \sup_{x:\ x\,\geq\,a_{k}(k\log k)^{1/2}}%
\ \bigg|\frac{\mathbf{P}(S_{k}\geq x)}{k\,\mathbf{P}(X_{1}\geq x)}%
\,-\,1\bigg|\,=\,0. \label{4.2}%
\end{equation}
Note that if $\,\delta>0,$\thinspace\ $k\geq\delta m^{n},$\thinspace\ and
$\,\varepsilon_{n}\geq\delta m^{-\varkappa n},$\thinspace\ then $\,\varepsilon
_{n}\geq\delta^{1+\varkappa}k^{-\varkappa}.$\thinspace\ Hence,%
\begin{equation}
\frac{\varepsilon_{n}k}{\left(  k\log k\right)  ^{1/2}}\ \geq\ \delta
^{1+\varkappa}\ \frac{k^{1/2-\varkappa}}{\left(  \log k\right)  ^{1/2}}\,.
\end{equation}
Since $\,0<\varkappa<1/2,$\thinspace\ the right hand side goes to infinity as
$\,k\uparrow\infty,$\thinspace\ and we will use it as $\,a_{k\,}.$%
\thinspace\ Thus, applying (\ref{4.2}) gives, as $\,n\uparrow\infty$,
\begin{gather}
\sum_{k>\delta m^{n}}\mathbf{P}(Z_{n}=k)\,\mathbf{P}(S_{k}\geq\varepsilon
_{n}k)\ =\ \left(  _{\!_{\!_{\,}}}1+o(1)\right)  \sum_{k>\delta m^{n}%
}k\,\mathbf{P}(Z_{n}=k)\,\mathbf{P}(X_{1}\geq\varepsilon_{n}k)\nonumber\\
\ =\ \left(  _{\!_{\!_{\,}}}1+o(1)\right)  a\,\varepsilon_{n}^{-\theta}%
\sum_{k>\delta m^{n}}k^{-(\theta-1)}\,\mathbf{P}(Z_{n}=k), \label{4.3}%
\end{gather}
where in the second step we used that $\,X_{1}$\thinspace\ has a tail of index
$\,\theta\in(2,1+\alpha).$ $\,$By (\ref{1.2}) we have
\[
\sum_{1\,\leq\,k\,\leq\,\delta m^{n}}k^{-(\theta-1)}\,\mathbf{P}%
(Z_{n}=k)\ \leq\ c\,m^{-\alpha n}\sum_{1\,\leq\,k\,\leq\,\delta m^{n}%
}k^{\alpha-\theta}\ \leq\ c\,m^{-(\theta-1)n}\,\delta^{1+\alpha-\theta}.
\]
By Theorem~1 of \cite{NeyVidyashankar2003}, for $\,\theta-1<\alpha,\,\ $we
have $\,\mathbf{E}\bigl\{Z_{n}^{-(\theta-1)};\,Z_{n}>0\bigr\}$ $\sim$
$I_{\theta}\,m^{-(\theta-1)n}$\thinspace\ as $n\uparrow\infty,$\thinspace
\ with $\,I_{\theta}$\thinspace\ defined in (\ref{notC.theta}). Hence, for all
sufficiently large $\,n,$
\begin{equation}
\Big|\sum_{k>\delta m^{n}}k^{-(\theta-1)}\,\mathbf{P}(Z_{n}=k)-I_{\theta
}\,m^{-(\theta-1)n}\Big|\ \leq\ c\,m^{-(\theta-1)n}\,\delta^{1+\alpha-\theta}.
\label{4.4}%
\end{equation}
Combining (\ref{4.3}) and (\ref{4.4}), we have the bound
\begin{equation}
\limsup_{n\uparrow\infty}\,\Big|\varepsilon_{n}^{\theta}\,m^{(\theta-1)n}%
\sum_{k>\delta m^{n}}\mathbf{P}(Z_{n}=k)\,\mathbf{P}(S_{k}\geq\varepsilon
_{n}k)-a\,I_{\theta}\Big|\ \leq\ c\,\delta^{1+\alpha-\theta}. \label{4.4'}%
\end{equation}
If $\,\varepsilon_{n}m^{\varkappa n}\rightarrow\infty$, then, obviously,
$\varepsilon_{n}^{-2\alpha}\,m^{-\alpha n}=o(\varepsilon_{n}^{-\theta
}\,m^{-(\theta-1)n})$. Therefore, by estimate (\ref{6.2}),
\begin{equation}
\limsup_{n\uparrow\infty}\,\varepsilon_{n}^{\theta}\,m^{(\theta-1)n}%
\sum_{1\,\leq\,k\,\leq\,\delta m^{n}}\mathbf{P}(Z_{n}=k)\,\mathbf{P}(S_{k}%
\geq\varepsilon_{n}k)\ \leq\ c\,\delta^{1+\alpha-\theta}. \label{4.1}%
\end{equation}
Part (b) follows from (\ref{4.4'}) and (\ref{4.1}) by letting $\,\delta
\downarrow0.$\smallskip

Finally, under $\,\varepsilon_{n}\sim\tau^{-1}m^{-\varkappa n},$%
\thinspace\ part (c) follows from (\ref{1.16}), (\ref{22}), (\ref{inL8}),
(\ref{6.1}), and (\ref{4.4'}). The proof is finished altogether.
\hfill$\square$

\subsection{Proof of Theorem~\ref{Bottcher}\label{SS.Bottcher}}

It follows from the assumed finiteness of an exponential moment of $X_{1}$,
see e.g. Lemma~III.5 in Petrov \cite{Petrov1975}, that for every $\delta
\in(0,1)$ there exists $h_{\delta}>0$ such that
\begin{equation}
\mathbf{E}\mathrm{e}^{hX_{1}}\,\leq\,\mathrm{e}^{\sigma^{2}(1+\delta)h^{2}%
/2},\qquad|h|\leq h_{\delta\,}.
\end{equation}
Thus, we may use the well-known Bernstein inequality, see Theorem~III.15 in
\cite{Petrov1975}. This gives, for all $k\geq1$ and $\varepsilon_{n}\leq
h_{\delta\,}$,
\begin{equation}
\mathbf{P}(S_{k}\geq\varepsilon_{n}k)\ \leq\ \exp\!\Big[-(1-\delta
)\,\frac{\varepsilon_{n}^{2}k}{2\sigma^{2}}\Big ].
\end{equation}
Therefore,
\begin{equation}
\mathbf{P}(R_{n}\geq\varepsilon_{n})\ \leq\ f_{n\!}\Bigl(\exp\!\Big[-(1-\delta
)\,\frac{\varepsilon_{n}^{2}}{2\sigma^{2}}\Big ]\Bigr)\quad\text{if
\ }\varepsilon_{n}\leq h_{\delta\,}. \label{5.1}%
\end{equation}
We may also assume that $\,\varepsilon_{n}\leq1/m.$\thinspace\ Set
$\,r_{n}:=\max\{k\geq1:\,m^{k}\leq\varepsilon_{n}^{-2}\}$. Then,
\begin{equation}
m^{-r_{n}-1}\,<\,\varepsilon_{n}^{2}\,\leq\,m^{-r_{n}}. \label{5.2}%
\end{equation}
The left hand inequality together with the monotonicity of $f_{n}$ gives
\begin{equation}
f_{n\!}\Bigl(\exp\!\Big[-(1-\delta)\,\frac{\varepsilon_{n}^{2}}{2\sigma^{2}%
}\Big ]\Bigr)\ \leq\ f_{n\!}\Bigl(\exp\!\Big[-(1-\delta)\,\frac{m^{-r_{n}-1}%
}{2\sigma^{2}}\Big ]\Bigr). \label{5.3}%
\end{equation}
Bounds (\ref{5.1}), (\ref{5.3}), and the right hand inequality in (\ref{5.2})
imply
\begin{equation}
\varepsilon_{n}^{-2\beta}m^{-n\beta}\log\mathbf{P}(R_{n}\geq\varepsilon
_{n})\ \leq\ \mu^{-n+r_{n}}\log f_{n\!}\Bigl(\exp\!\Big[-(1-\delta
)\,\frac{m^{-r_{n}-1}}{2\sigma^{2}}\Big ]\Bigr), \label{5.3'}%
\end{equation}
where we used $\mu=m^{\beta}.$\thinspace\ Since $\,r_{n}\rightarrow\infty
,$\thinspace\ by the Kesten-Stigum theorem for supercritical Galton-Watson
processes,
\begin{equation}
\lim_{n\uparrow\infty}f_{r_{n}+1}\Big(\exp\!\Big[-(1-\delta)\,\frac
{m^{-r_{n}-1}}{2\sigma^{2}}\Big ]\Bigr)\ =\ \varphi\bigl((1-\delta
)/2\sigma^{2}\bigr). \label{5.4}%
\end{equation}
On the other hand, from the assumption $\,\varepsilon_{n}^{2}m^{n}%
\rightarrow\infty$\thinspace\ and the right hand inequality in (\ref{5.2}) it
follows that $\,n-r_{n}\rightarrow\infty.$\thinspace\ Therefore, by
(\ref{BConv}) we have for $s\in\lbrack0,1]$,
\begin{equation}
\lim_{n\uparrow\infty}\mu^{-n+r_{n}+1}\log f_{n-r_{n}-1}(s)\ =\ \log
\mathsf{B}(s). \label{5.5}%
\end{equation}
By the continuity of $\,\mathsf{B},$\thinspace\ combining (\ref{5.4}) and
(\ref{5.5}) we obtain
\begin{equation}
\lim_{n\uparrow\infty}\mu^{-n+r_{n}+1}\log f_{n\!}\Bigl(\exp\!\Big[-(1-\delta
)\,\frac{m^{-r_{n}-1}}{2\sigma^{2}}\Big ]\Bigr)\ =\ \log\mathsf{B}%
\bigl(\varphi\bigl((1-\delta)/2\sigma^{2}\bigr)\bigr). \label{5.6}%
\end{equation}
Now (\ref{UB'}) follows from (\ref{5.3'}) and (\ref{5.6}) letting
$\,\delta\downarrow0$.\smallskip

In order to prove (\ref{LB'}) we will exploit the following version of
Kolmogorov's inequality: For $\,0<\delta<1$\thinspace\ fixed, there exists a
constant $D\in(0,\infty)$ such that
\begin{equation}
\mathbf{P}(S_{k}\geq\varepsilon_{n}k)\ \geq\ \exp\!\Big[-(1+\delta
)\,\frac{\varepsilon_{n}^{2}k}{2\sigma^{2}}\Big ],\qquad k>D/\varepsilon
_{n\,}^{2},\quad n\geq1. \label{5.7}%
\end{equation}
See Statulevicius \cite{Statulevicius1966}. Using (\ref{5.7}) we obtain%
\begin{gather}
\mathbf{P}(R_{n}\geq\varepsilon_{n})\ \geq\ \sum_{k>D/\varepsilon_{n}^{2}%
}\mathbf{P}(Z_{n}=k)\,\exp\!\Big[-(1+\delta)\,\frac{\varepsilon_{n}^{2}%
k}{2\sigma^{2}}\Big ]\nonumber\\
\geq\ f_{n\!}\!\left(  \exp\!\Big[-(1+\delta)\,\frac{\varepsilon_{n}^{2}%
}{2\sigma^{2}}\Big ]\right)  -\,\mathbf{P}(Z_{n}\leq D/\varepsilon_{n}^{2}).
\label{5.7''}%
\end{gather}
Clearly, if $\,D/\varepsilon_{n}^{2}<\mu^{n},$\thinspace\ then $\,\mathbf{P}%
(Z_{n}\leq D/\varepsilon_{n}^{2})=0,$\thinspace\ and we pass directly to
statement (\ref{5.10}) below. Otherwise, it follows from
Proposition~\ref{P.Boettcher} that
\begin{equation}
\mathbf{P}(Z_{n}\leq D/\varepsilon_{n}^{2})\ \leq\ \exp\!\Big[-c\,D^{-\beta
/(1-\beta)}\,(\varepsilon_{n}^{2}m^{n})^{\beta/(1-\beta)}\Big ]. \label{5.7'}%
\end{equation}
\textrm{F}rom (\ref{5.7''}), (\ref{5.7'}), and the left hand inequality in
(\ref{5.2}), we have
\begin{equation}
\mathbf{P}(R_{n}\geq\varepsilon_{n})\ \geq\ f_{n\!}\Bigl(\exp\!\Big[-(1+\delta
)\,\frac{m^{-r_{n}}}{2\sigma^{2}}\Big ]\Bigr)-\exp\!\Big[-c\,(\varepsilon
_{n}^{2}m^{n})^{\beta/(1-\beta)}\Big ]. \label{5.8}%
\end{equation}
Analogously to (\ref{5.6}),
\begin{equation}
\lim_{n\uparrow\infty}\mu^{-n+r_{n}}\log f_{n\!}\Bigl(\exp\!\Big[-(1+\delta
)\,\frac{m^{-r_{n}}}{2\sigma^{2}}\Big ]\Bigr)\ =\ \log\mathsf{B}%
\bigl(\varphi\bigl((1+\delta)/2\sigma^{2}\bigr)\bigr). \label{5.9}%
\end{equation}
By the left hand inequality of (\ref{5.2}), $\,\mu^{n-r_{n}}\leq m^{\beta
}(\varepsilon_{n}^{2}m^{n})^{\beta}.$\thinspace\ Therefore, from the limit
statement (\ref{5.9}) we see that the second term at the right hand side of
estimate (\ref{5.8}) is negligible compared with the first term there, i.e.
\begin{equation}
\mathbf{P}(R_{n}\geq\varepsilon_{n})\ \geq\ f_{n\!}\Bigl(\exp\!\Big[-(1+\delta
)\,\frac{m^{-r_{n}}}{2\sigma^{2}}\Big ]\Bigr)(1+o(1)). \label{5.10}%
\end{equation}
Thus, using the left hand inequality in (\ref{5.2}), we get the bound
\begin{equation}
\varepsilon_{n}^{-2\beta}m^{-n\beta}\log\mathbf{P}(R_{n}\geq\varepsilon
_{n})\,\geq\,\mu^{-n+r_{n}+1}\log f_{n\!}\Bigl(\exp\!\Big[-(1+\delta
)\,\frac{m^{-r_{n}}}{2\sigma^{2}}\Big ]\Bigr)+\,o(1). \label{5.8'}%
\end{equation}
Since $\delta$ is arbitrary, combining (\ref{5.8'}) and (\ref{5.9}) completes
the proof of (\ref{LB'}).\smallskip

In the derivation of (\ref{5.10}) from (\ref{5.7''}) we learned that the
second term at the right hand side of (\ref{5.7''}) is small compared with the
first term there. Thus, from (\ref{5.7''}) together with (\ref{5.1}) we get
\begin{align}
f_{n\!}\Bigl(\exp\!\Big[-(1+\delta)\,\frac{\varepsilon_{n}^{2}}{2\sigma^{2}%
}\Big ]\Bigr)\left(  _{\!_{\!_{\,}}}1+o(1)\right)  \  &  \leq\ \mathbf{P}%
(R_{n}\geq\varepsilon_{n})\\
&  \leq\ f_{n\!}\Bigl(\exp\!\Big[-(1-\delta)\,\frac{\varepsilon_{n}^{2}%
}{2\sigma^{2}}\Big ]\Bigr).\nonumber
\end{align}
Hence, if $\,\varepsilon_{n}^{2}=m^{-\lambda_{n}}$\thinspace\ then
(\ref{BottAsymp}) follows from these inequalities and (\ref{5.9}) replacing
there $\,r_{n}\,\ $by $\,\lambda_{n\,},$\thinspace\ and finally letting
$\,\delta\downarrow0.$\thinspace\ Altogether, the proof of
Theorem~\ref{Bottcher} is complete.\hfill$\square$

\subsection{Proof of Theorem~\ref{T_LDev1}}

With $\,B_{2}$\thinspace\ from Proposition~\ref{P.Boettcher}, and $\,\theta
>2$\thinspace\ the tail index of $\,X_{1\,},$\thinspace\ define $\,k_{n}%
:=m^{n}/\log^{(1-\beta)/\beta}m^{2n\theta/B_{2}}.$\thinspace\ Then by
Proposition~\ref{P.Boettcher}, for all sufficiently large $\,n,$
\begin{equation}
\mathbf{P}(Z_{n}\leq k_{n})\ \leq\ \exp\bigl[-(B_{2}/2)(k_{n}/m^{n}%
)^{-\beta/(1-\beta)}\bigr]\ =\ m^{-\theta n}.
\end{equation}
Hence, for these $\,n,$
\begin{equation}
\sum_{k\leq k_{n}}\mathbf{P}(Z_{n}=k)\,\mathbf{P}(S_{k}\geq\varepsilon
_{n}k)\ \leq\ \mathbf{P}(Z_{n}\leq k_{n})\ \leq\ m^{-\theta n} \label{4.5}%
\end{equation}
and
\begin{equation}
\sum_{k\leq k_{n}}k^{-(\theta-1)}\,\mathbf{P}(Z_{n}=k)\ \leq\ \mathbf{P}%
(Z_{n}\leq k_{n})\ \leq\ m^{-\theta n}. \label{4.6}%
\end{equation}
It is easy to verify that
\begin{equation}
\frac{\varepsilon_{n}k_{n}}{(k_{n}\log k_{n})^{1/2}}\ =\ \left(
_{\!_{\!_{\,}}}c+o(1)\right)  \varepsilon_{n}m^{n/2}n^{-1/2\beta}\quad\text{as
}\,n\uparrow\infty.
\end{equation}
By our assumption in the theorem, the right hand side converges to infinity.
Then, we can use (\ref{4.2}) with $\,a_{k}:=\varepsilon_{n}(k/\log k)^{1/2}%
$\thinspace\ to obtain
\begin{gather}
\sum_{k>k_{n}}\mathbf{P}(Z_{n}=k)\,\mathbf{P}(S_{k}\geq\varepsilon
_{n}k)=\ \left(  _{\!_{\!_{\,}}}1+o(1)\right)  \sum_{k>k_{n}}k\,\mathbf{P}%
(Z_{n}=k)\,\mathbf{P}(X_{1}\geq\varepsilon_{n}k)\nonumber\\
\ =\ \left(  _{\!_{\!_{\,}}}1+o(1)\right)  a\,\varepsilon_{n}^{-\theta}%
\sum_{k>k_{n}}k^{-(\theta-1)}\,\mathbf{P}(Z_{n}=k)\quad\text{as }%
\,n\uparrow\infty. \label{4.7}%
\end{gather}
Theorem~1 of \cite{NeyVidyashankar2003} and (\ref{4.6}) yield
\begin{equation}
\sum_{k>k_{n}}k^{-(\theta-1)}\,\mathbf{P}(Z_{n}=k)\ =\ I_{\theta}%
\,m^{-(\theta-1)n}\left(  _{\!_{\!_{\,}}}1+o(1)\right)  \quad\text{as
}\,n\uparrow\infty.
\end{equation}
Substituting this into (\ref{4.7}) and combining with (\ref{4.5}) completes
the proof. \hfill$\square$\bigskip

{\small
\bibliographystyle{alpha}
\bibliography{bibtex,bibtexmy}
}


\vfill
\end{document}